\documentclass{article}
\usepackage[utf8]{inputenc}
\usepackage[utf8]{inputenc}
\usepackage{amsmath}
\usepackage[left=3.0cm, right=3.0cm, bottom=2.5cm, top=2.5cm]{geometry}
\usepackage{hyperref}
\usepackage{cleveref}
\usepackage{amssymb}
\usepackage{mathtools}
\usepackage{tikz-cd}
\usepackage{amsthm}
\usepackage{graphicx}
\usepackage{graphics}
\usepackage{mathrsfs}
\usepackage{cleveref}
\usepackage{thmtools}
\usepackage{enumitem}
\usepackage{bbm}
\newlist{steps}{enumerate}{1}

\usepackage{array}
\setlist[steps, 1]{label = Step \arabic*:}
\theoremstyle{plain}
\usepackage{footnote}
\makesavenoteenv{tabular}
\newtheorem*{theorem*}{Theorem}
\newtheorem{theorem}{Theorem}[subsection]
\newtheorem{definition}[theorem]{Definition}

\newtheorem{remark}[theorem]{Remark}

\newtheorem{lemma}[theorem]{Lemma}

\numberwithin{equation}{subsection}
\usepackage{rotating}
\usepackage{comment}
\usepackage[nottoc,notlot,notlof]{tocbibind} 

\usepackage[toc,page]{appendix}

\title{A Canonical Quantization of Poisson Manifolds: a 2-Groupoid Scheme} 
\author{Joshua Lackman\footnote{josh@pku.edu.cn} }
\date{}
\begin{document}

\maketitle
\begin{abstract}
\noindent We canonically quantize a Poisson manifold to a Lie 2-groupoid, complete with a quantization map, and show that it relates geometric and deformation quantization: the perturbative expansion in $\hbar$ of the (formal) convolution of two quantized functions yields Kontsevich's star product. Meanwhile, we can push forward this quantization map (by integrating over homotopies of paths) to obtain a quantization map in traditional geometric quantization. This gives a polarization-free, path integral definition of the quantization map, which does not have a prior prescription for Poisson manifolds and which only has a partial prescription for symplectic manifolds. We show that in conventional quantum mechanics our operators are defined on the entire prequantum Hilbert space and that they all preserve polarized sections, contrary to the Kostant-Souriau approach.
\end{abstract}
\tableofcontents
\section{Introduction}
 Kontsevich proved that every Poisson manifold $(M,\Pi)$ admits a formal deformation quantization,\footnote{Obtained using a sequence of bidifferential operators.} called a star product (see \cite{kontsevich}). That is, there is an associative product $\star$ on $C^{\infty}(M)[[\hbar]]$ which deforms the product on $C^{\infty}(M),$ in such a way that 
\begin{equation}\label{star}
    f\star g=fg+i\hbar\{f,g\}+\mathcal{O}(\hbar^2)\;.
    \end{equation}
This is a perturbative approach to quantization. We can think of the inclusion map
\begin{equation}
    C^{\infty}(M)\xhookrightarrow{} C^{\infty}(M)[[\hbar]]
\end{equation}
as a perturbative \textit{quantization map}, which associates to $f$ the operator $g\mapsto f\star g.$ Since this is a perturbative approach we don't get a $C^*$-algebra, which is important in physics.
\\\\A non-perturbative approach to quantization, which was initiated by Weinstein in \cite{weinstein}, is geometric quantization of Poisson manifolds. This quantization method aims to define a $C^*$-algebra of observables — it is given by polarized sections of a twisted convolution algebra of the symplectic groupoid (the groupoid which integrates the Poisson structure, \cite{eli}). In the symplectic case, we would like this $C^*$-algebra to act on polarized sections of the prequantum line bundle.\footnote{For a discussion on when a star product acts on the space of sections of such a quantization, see \cite{nolle}, \cite{nolle2}.}
\\\\Geometric quantization is much less understood than formal deformation quantization (which is consistent with the current understanding of QFT) and is incomplete as it doesn't prescribe a quantization map.\footnote{Here, we mean an injective $^*$-linear map from $C^{\infty}_c(M)$ to the $C^*$-algebra of observables that determines a star product.} In this sense, the difficulty of geometric quantization is complementary to the difficulty of formal deformation quantization: one (partially) prescribes a $C^*$-algebra but not a quantization map, the other prescribes a formal quantization map but not a $C^*$-algebra. 
\\\\In some examples (eg. \cite{pierre},\cite{bone},\cite{reiffelact}) it has been shown that there exists a quantization map which is compatible with formal deformation quantization, but an explanation of why and when this occurs is lacking. Even in traditional geometric quantization of symplectic manifolds the quantization map is, at best, defined only on a small subspace (and assumes a polarization exists). It is desirable to have a more mathematically precise theory which directly relates these two quantization methods.
\\\\The constructions we make are justified by results of the Poisson sigma model in \cite{bon}, \cite{catt}. We \textit{formally} define a quantization map (\cref{qant}) by using the Lie 2-groupoid $\Pi_2(T^*M)$ integrating a Poisson manifold (\cite{zhuc}) — this involves computing a path integral over homotopies of algebroid paths. In the symplectic case, this 2-groupoid is a quotient of the singular simplicial set of the manifold. Noncommutativty for $\hbar\ne 0$ arises due to distinct points in phase space being isomorphic in the 2-groupoid. This is a two dimensional sigma model explanation of quantization, somewhat as in \cite{guk}.
\subsection{Main Example}
We compute the path integrals in the case of the pointed K\"{a}hler manifold $(\textup{T}^*\mathbb{R},\omega,J,0)$ with the canonical symplectic form (\cref{exam}).\footnote{The only examples of geometric quantization of symplectic manifolds that have been shown to have an injective quantization map are discrete quotients of $\textup{T}^*\mathbb{R}^n,$ and it is not difficult to compute the path integrals there.} The quantization of a smooth function $f$ is a section of a trivializable line bundle over the symplectic groupoid, and can be identified with (\cref{geoq})
\begin{align}
     Q_f:\textup{T}^*\mathbb{R}\times \textup{T}^*\mathbb{R}\to\mathbb{C}\;,\,\; Q_f(u,v)=\frac{1}{4\pi\hbar}\int_{\textup{T}^*\mathbb{R}}f(z)e^{\frac{i}{\hbar}P_0(u,v,z)}\,\omega\;,
\end{align}
where $P_0(u,v,z)$ is the (signed) area of the polygon with vertices $0,u,v,z \in \textup{T}^*\mathbb{R}\,,$ and the integration is over $z.$
\\\\The canonical connection on the prequantum line bundle is $(pdq-qdp)/2$ (\cref{multi0}) and parallel transport determines a representation of these quantized functions as integral operators on $L^2(\textup{T}^*\mathbb{R}).$ The result is 
\begin{align}\label{defi}
    Q_f\Psi(u)=\frac{1}{(4\pi\hbar)^2}\int_{\textup{T}^*\mathbb{R}\times \textup{T}^*\mathbb{R}}f(v)\Psi(z)e^{\frac{i}{\hbar}\Omega(u,v,z)}\,\omega\boxtimes\omega\;,
\end{align}
where $\Psi\in L^2(\textup{T}^*\mathbb{R})$ and $\Omega(u,v,z)$ is the (signed) area of the geodesic triangle determined by $u,v,z;$ the integral is over $v,z,$ with respect to the product measure.
\\\\The prequantum Hilbert space $L^2(\textup{T}^*\mathbb{R})$ contains both pure and mixed states.\footnote{In the $C^*$-algebra sense. A pure vector state is one which generates an irreducible representation (\cite{Gleason}).} Polarized sections form irreducible subrepresentations for \textit{all} covariantly constant Lagrangian polarizations. That is, if $\Psi$ is polarized (with respect to some complex linear polarization), then (\cref{ch})
\begin{equation}
Q_f\Psi(u)=\frac{1}{4\pi\hbar}\int_{\textup{T}^*\mathbb{R}}f(u')\Psi(u'-u)e^{\frac{i}{\hbar}\Omega(0,u,u')}\,\omega\;,
\end{equation}
and this is a polarized section. In particular, we get the position, momentum and Segal-Bargmann (holomorphic) representations. The states of the latter take the form 
\begin{equation}
    \Psi(p,q)=e^{-\frac{|p+iq|^2}{2\hbar}}\psi(p+iq)\;.
\end{equation}
With the Kostant–Souriau prescription, only a small subspace of functions prequantize to operators which preserve polarized sections, and this subspace depends on the polarization  — \textit{only} affine functions prequantize to operators which preserve all linearly polarized sections. Our operators \ref{defi} are defined \textit{polarization-free} and they \textit{all} preserve polarized sections, for \textit{all} complex linear polarizations. These two quantizations agree on the subspace where the Kostant–Souriau prescription is defined.
\\\\In addition, with respect to the twisted convolution algebra $*$ on the groupoid and the non-perturbative Moyal product $\star$ (\cref{bon}, \cref{my}), $Q_{f\star g}\Psi=(Q_f\ast Q_g)\Psi=Q_f(Q_g\Psi).$\footnote{The image of $\Omega$ under the van Est map is $\omega$ (\cite{Lackman4}) and is therefore independent of the complex structure.}
\begin{remark}
We view the role of Lagrangian polarizations in the usual definitions of geometric quantization as a tool to help guess a quantization map, but not necessarily as an essential ingredient (sometimes they don't even exist). This is consistent with Kontsevich's and Fedosov's (\cite{fedosov}) constructions of the star product, where polarizations are absent. We will only use them to help construct irreducible representations.
\\\\In the usual approach, a common heuristic is that the prequantum Hilbert space being ``too large" is a problem. However, our perspective is that the only problem is that the subalgebra of operators generated by the Kostant-Souriau prequantization map is far too large to represent observables — its image only contains first order differential operators. A precise version of this is that the Kostant-Souriau prequantization map is not compatible with a star product (\cref{nd}). This perspective is particularly useful when the prequantum Hilbert space contains inequivalent irreducible representations.
\end{remark}
\subsection{Lie 2-Groupoid Approach to Quantization and Main Results}\label{main}
Our 2-groupoid approach to quantization relates formal deformation quantization and geometric quantization of Poisson manifolds, and is closely related to the Poisson sigma model. It has the advantage that it can be defined on a lattice. The idea is based on the following: every Poisson manifold $(M,\Pi)$ has a canonical integration given by a Lie 2-groupoid, denoted $\Pi_2(\textup{T}^*M)$ (see \cite{zhuc}). This Lie 2-groupoid comes with a canonical multiplicative prequantum line bundle 
\begin{equation}
    \mathcal{L}_2\to \Pi_2(\textup{T}^*M)\;.
\end{equation}
Its space of sections, denoted $\Gamma(\mathcal{L}_2),$ is formally a $C^*$-algebra — its product $*$ is given by convolution (this involves a path integral) and its involution $^*$ is complex conjugation. It's norm is (formally) the $L^2$-norm of the left regular representation.\footnote{Every algebra acts on itself on the left, forming a representation called the left regular representation.}
There is well-defined, canonical quantization map
\begin{equation}
    q:C^{\infty}_c(M)\to \Gamma(\mathcal{L}_2)\;.
\end{equation} 
In this paper we regard a Poisson manifold as a smooth family of classical systems, and therefore we want a corresponding smooth family of quantizations (which often exist formally, but not always. See \cref{non}). In this context, we take the pair 
\begin{equation}
    (\mathcal{L}_2\to\Pi_2(\textup{T}^*M),\,q)
    \end{equation}
to be a ``universal" quantization of $(M,\Pi).$ We have the following theorem:
\begin{remark}
We want to emphasize that our use of the Lie 2-groupoid has nothing to do with the existence of the symplectic groupoid — we need to use the 2-groupoid in order to define $q.$
\end{remark}
\begin{theorem}(formal)\label{formal} Let $(M,\Pi)$ be a Poisson maniold such that the convolution of two sections of $\Gamma(\mathcal{L}_2)$ is smooth.\footnote{This should be true in a wide class of examples, including regular Poisson manifolds satisfying an integrality condition.} We have that
\begin{enumerate}
    \item For $m\in M,$ $(q_f\ast q_g)(m)$ is equal to the 3-point function
    \begin{equation}
 \int_{X(\infty)=m} f(X(1))g(X(0))\,e^{\frac{i}{h}S[X,\eta]}\,\mathcal{D}X\,\mathcal{D}\eta
\end{equation}of the Poisson sigma model on the disk (\cite{catt}). Therefore, perturbatively at $\hbar=0,$ 
    \begin{equation}
        (q_f\ast q_g)\vert_M\sim f\star g.
        \end{equation}
    \item In the case that the symplectic groupoid $\Pi_1(\textup{T}^*M)$ is prequantizable, $\mathcal{L}_2\to \Pi_2(\textup{T}^*M)$ descends to a multiplicative prequantum line bundle 
    \begin{equation}
        \mathcal{L}_1\to\Pi_1(\textup{T}^*M)
        \end{equation}
        and $q$ pushes forward to a ``geometric" quantization map
    \begin{equation}
        Q:C^{\infty}_c(M)\to \Gamma(\mathcal{L}_1)\;.
    \end{equation}
    \item In the case that $(M,\Pi)$ is pointed, simply connected and symplectic, there is a canonical representation of $\Gamma(\mathcal{L}_1)$ on sections of the prequantum line bundle $\mathcal{L}_0\to M,$ hence there is a canonical representation of the $C^*$-algebra generated by $Q.$ 
\end{enumerate}
\end{theorem}
$\,$\\To clarify the different (but related) quantization maps:
\begin{enumerate}
    \item The map $q$ is the 2-groupoid quantization map. The others are obtained from it.
    \item The map $Q$ in part 2 is the geometric quantization map.
    \item In part 1, the left regular representation defines a third quantization map $\mathcal{Q},$  
\begin{equation}\label{quant}
    f\mapsto\mathcal{Q}_f\,,\;\;\mathcal{Q}_fg=(q_f\ast q_g)\vert_M\;.
\end{equation} 
We describe the $C^*$-algebra which $\mathcal{Q}_f$ belongs to in the next section, see \cref{product}; we denote it by  $\mathcal{B}(\Pi).$ 
\item From the representation in part 3, we get a fourth quantization map by composing $Q$ with the map into the bounded linear operators on the Hilbert space of the prequantum line  bundle.
\end{enumerate}
The main goal is to realize \cref{formal} non-perturbatively, which we are optimistic about in the case of regular Poisson manifolds. This is the goal of the lattice framework recently developed in \cite{Lackman4} for defining path integrals over spaces of Lie algebroid morphisms, see appendix \ref{appen}. Most known examples of non-perturbative deformation quantizations fit into this framework, as discussed in example 9.2.7 of the aforementioned paper. In conventional quantum mechanics, we work out this result in \cref{exam}.
\\\\In the next section we will make clear exactly what we mean by a \textit{quantization}.
\subsection{Non-Perturbative (Tangential) Deformation Quantizations}\label{non}
We will discuss non-perturbative deformation quantizations. The definition we use is essentially the most general one that is compatible with both $C^*$-algebraic quantization and formal deformation quantization. It sits between Rieffel's strict quantization (\cite{rieffel}) and his strict deformation quantization.
\begin{definition}\label{nd}(see \cite{eli2})
Let $(M,\Pi)$ be a Poisson manifold and let $A\subset[0,1]$ be a set containing $0$ as an accumulation point. For each $\hbar\in A,$ let $M_{\hbar}$ be a unital $C^*$-algebra such that $M_{0}=L^{\infty}(M)$ and let 
\begin{equation}
    Q_{\hbar}:C_c^{\infty}(M)\to M_{\hbar}
\end{equation}
be injective and $^*$- linear such that its image generates $M_{\hbar}\,.$ Furthermore, assume that $Q_0$ is the inclusion map. We say that $Q_{\hbar}$ is a (non-perturbative) deformation quantization of $(M,\Pi)$ if there is a star product $\star_{\hbar}$ on $C^{\infty}(M)[[\hbar]]$ such that, for all $n\in\mathbb{N},$
\begin{equation}
    \frac{1}{\hbar^{n}}||Q_{\hbar}(f)Q_{\hbar}(g)-Q_{\hbar}(f\star_{\hbar}^n g)||_{\hbar}\xrightarrow[]{\hbar \to 0} 0\;,
\end{equation}
where $f\star_{\hbar}^n g$ is the component of the star product up to order $n.$
\end{definition}
$\,$\\In this paper we view a Poisson manifold as a smooth family of classical systems, thus we want a quantization which quantizes each of the classical systems. This motivates the next definition:
\begin{definition}
We say that a (non-perturbative) deformation quantization of $(M,\Pi)$ is tangential if for all $\hbar\in A,$ $Q_{\hbar}(f)Q_{\hbar}(g)=0$ whenever $f,g$ are supported on disjoint symplectic leaves.
\end{definition}
Every regular Poisson manifold has a tangential, formal deformation quantization,\footnote{A star product is tangential if the sequence of bidifferential operators are tangent to the symplectic leaves.} see the work of Fedosov \cite{fedosov}, also see  \cite{gamme}, \cite{wein}. (there are a number of examples of non-regular Poisson manifolds which also have such a quantization, as discussed in \cite{Lackman4}).\footnote{A quantization which isn't tangential exhibits a kind of tunneling.} The maps discussed in \cref{main} have a tangential property since the orbits of the Lie 2-groupoid $\Pi_2(\textup{T}^*M)$ are the symplectic leaves.
\\\\A Poisson manifold has a natural $C^*$-algebra associated to it, which in the case of a symplectic manifold is the space of bounded linear operators on square integrable functions. For a Poisson manifold, it is essentially the direct product of the $C^*$-algebras of its symplectic leaves:
\begin{definition}\label{product}
Let $(M,\Pi)$ be a Poisson manifold. We have a norm on the subspace of measurable functions for which
\begin{equation}
    ||f||:=\textup{sup}_{\mathcal{L}}||f||_{\mathcal{L}}<\infty\;,
\end{equation}
where the supremum is over all symplectic leaves $\mathcal{L}$ and $||f||_{\mathcal{L}}$ is the $L^2$-norm of $f$ pulled back to $\mathcal{L}$\footnote{With respect to the measure induced by the symplectic form.} We denote this space by $L^2(\Pi).$
\\\\Let $\mathcal{B}(\Pi)$ denote the space of bounded linear operators $T$ on $L^2(\Pi)$ such that, if $f\vert_{\mathcal{L}}=0$ for a symplectic leaf $\mathcal{L},$
then \begin{equation}
    (Tf)\vert_{\mathcal{L}}=0\;.
\end{equation}There is an involution on $\mathcal{B}(\Pi)$ defined by
\begin{equation}
    \textup{T}^*\vert_{\mathcal{L}}=\big(T\vert_{\mathcal{L}}\big)^*\;.
    \end{equation}
This turns $\mathcal{B}(\Pi))$ into a $C^*$-algebra. 
\end{definition}
$\,$\\If $\Pi=0,$ then $L^2(\Pi)=L^{\infty}(M).$ 
\\\\Motivated by \cref{formal} and \cref{quant}, we are optimistic that regular Poisson manifolds such that the Poisson structure is integral over spheres (ie. satisfies \ref{cond}) have non-perturbative, tangential deformation quantizations for which 
\begin{equation}
  M_{\hbar}\subset \mathcal{B}(\Pi)\;.
\end{equation}
This would give a a non-perturbative form of Kontsevich's star product on such Poisson manifolds.\footnote{For examples of non-tangential quantizations, see the quantizations of the dual of nilpotent Lie algebras \cite{rieffel2}.}
\begin{remark}
In formal deformation quantization, it's really the gauge equivalence class of the formal quantization map that's important. Similarly, we might define two quantization maps to be gauge equivalent if there is an automorphism of $C^*$-algebras intertwining them. 
\end{remark}
\section{The Star Product and the Poisson Sigma Model}
A Poisson sigma model construction of Kontsevich's star product as a 3-point function in a QFT was explicitly written out by Cattaneo and Felder in \cite{catt}.
\\\\Precisely, the domain of this Poisson sigma model is the disk $D$ with three cyclically ordered marked points on the boundary, denoted $0\,,1\,,\infty.$ The target space is a Poisson manifold $(M,\Pi)$ and a field is a vector bundle morphism 
\begin{equation*}
    (X,\eta):TD\to \textup{T}^*M
\end{equation*}
covering $X:D\to M\,.$ There is a boundary condiition which says that the pullback of $\eta$ to $\partial D$ takes values in the zero section. The action is given by
\vspace{0.05cm}\\
\begin{equation}\label{action}
    S[X,\eta]=\int_D \eta\wedge dX+\frac{1}{2}\Pi\vert_X(\eta\,,\eta)\;.
\end{equation} The star product is then given by\footnote{More precisely, we can say that the perturbative expansions at $\hbar=0$ is Kontsevich's star product.}
\vspace{0.05cm}\\
\begin{equation}\label{star}
 (f\star g)(m)=\int_{X(\infty)=m} f(X(1))g(X(0))\,e^{\frac{i}{h}S[X,\eta]}\,\mathcal{D}X\,\mathcal{D}\eta\;.
\end{equation}
\vspace{0.05cm}\\
Here, $f,g\in C^{\infty}(M)$ and the integral is over all fields $X,\eta$ with $X(\infty)=m.$ The measure is normalized so that $1\star 1=1.$
\subsection{The Tangential Star Product}
Here, we discuss the tangential property of \cref{star}. In the following lemma, $\pi:\textup{T}^*M\to M$ is the projection and the domain of integration is the set of Lie algebroid morphisms $X:TD\to \textup{T}^*M\,,$ where the Lie algebroid structure on $ \textup{T}^*M$ is determined by the Poisson structure $\Pi,$ ie. the bracket on exact sections of $\textup{T}^*M$ is 
$[df,dg]=d\{f,g\},$ and for $\alpha\in\textup{T}^*M$ the anchor is $\alpha\mapsto \Pi(\alpha,\cdot)$ (\cite{ruif}). 
\\\\The action is given by 
\begin{equation}
S[X]=\int_D X^*\Pi\;.
\end{equation}
\begin{lemma}\label{bon}(see \cite{bon})
In the following, if the quantity on the right is smooth on $M,$ then we have that
\begin{equation}\label{regstar}
    (f\star g)(m)=\int_{\pi(X(\infty))=m}f(\pi(X(1)))\,g(\pi(X(0)))\,e^{\frac{i}{\hbar}S[X]}\,\mathcal{D}X\;.
    \end{equation}
    \end{lemma}
The right side of \ref{regstar} is exactly what one obtains by restricting the domain of integration of \ref{star} to Lie algebroid morphisms, ie. classical solutions of the Euler-Lagrange equations. This is tangential because Lie algebroid morphisms stay in a symplectic leaf.
\\\\The problem of computing $f\star g$ non-perturbatively is to compute (as a distribution on $M\times M$) the integral kernel 
\begin{equation}\label{kernel}
  (m_0,m_1,m)\mapsto  \int_{m_0,m_1,m} e^{\frac{i}{\hbar}S[X]}\,\mathcal{D}X\;,
\end{equation}
where the domain of integration is the set of Lie algebroid morphisms $X$ such that
\begin{equation}
    \pi(X(0))=m_0,\,\pi(X(1))=m_1,\, \pi(X(\infty))=m\;.
\end{equation}
A lattice approach to defining these path integrals is described in \cite{Lackman4}, see appendix \ref{appen}. On a K\"{a}hler manifold this kernel should be canonically defined, as there is a canonical 2-cocycle on the local pair groupoid (this is consistent with the fact that Fedosov's quantization is canonically defined on a K\"{a}hler manifold, as the Levi-Civita connection is a symplectic connection). See \cref{moyal}.
\section{ Geometric Quantization of the Lie 2-Groupoid}
Based on \cref{regstar}, we are going to interpret \cref{star} using the twisted convolution algebra of a Lie 2-groupoid, obtained via a geometric quantization-like construction. Much of this was explained in \cite{Lackman3}.
\\\\If $(M,\Pi)$ is a Poisson manifold then $\textup{T}^*M$ is naturally a Lie algebroid, as discussed in the previous section. We can integrate it to a Lie 2-groupoid.
\begin{definition}
We let $|\Delta^n|$ denote the geometric realization of the standard n-simplex.
\end{definition}
We now define the 2-groupoid, see \cref{def} for more details
\begin{definition}(see \cite{zhuc})\label{lie2}
Let $(M,\Pi)$ be a Poisson manifold. There is a Lie 2-groupoid, denoted $\Pi_{2}(\textup{T}^*M),$ which in degree $0$ is given by $M,$ in degree $1$ is given by
\begin{equation}
  \Pi_{2}^{(1)}(\textup{T}^*M)=\{\gamma\in \textup{Hom}(\,\textup{T}[0,1],\, \textup{T}^*M\,): \gamma\vert_{\{0,1\}}=0\}\;,
\end{equation}
and in degree $2$ is given by 
\begin{equation}
\Pi_{2}^{(2)}(\textup{T}^*M)=\{X\in \textup{Hom}(\,\textup{T}|\Delta^2|,\, \textup{T}^*M\,): X\vert_{\textup{vertices}}=0\}/\sim \;,
\end{equation} where $\sim$ identifies two morphisms if they are homotopic relative to the boundary. We denote the face maps 
\begin{equation}
    \Pi_{2}^{(1)}(\textup{T}^*M)\to M
\end{equation}
by $s,t,$ where $s(\gamma)=\gamma(0),\,t(\gamma)=\gamma(1).$ We denote the
face maps 
\begin{equation}
   \Pi_{2}^{(2)}(\textup{T}^*M)\to \Pi_{2}^{(1)}(\textup{T}^*M)
\end{equation}
by $ d_0,\,d_1,\,d_2,$ where 
\begin{equation}
    d_0(X)=\gamma_1,\;d_1(X)=\gamma_2,\;d_2(X)=\gamma_0
    \end{equation}
are the 1-dimensional faces of $X.$ For such an $X,$ we may write 
\begin{equation}
    \gamma_0\cdot_X\gamma_1=\gamma_2.
    \end{equation}
\end{definition}$\,$
\begin{remark}\label{def}
    In this definition, $X\vert_{\textup{vertices}}=0$ means that $X$ takes values in the zero section when restricted to the pullback of the tangent bundle to the vertices. Lie algebroid morphisms $X_0,X_1$ are homotopic relative to the boundary if there is a Lie algebroid morphism 
\begin{equation}
    T|\Delta^2|\times T[0,1]\to \textup{T}^*M
    \end{equation}
which agrees with $X_0, X_1$ on  $T|\Delta^2|\times \{0\}\,,\,T|\Delta^2|\times \{1\},$ respectively, and which takes values in the zero section on \begin{equation}
    T\partial|\Delta^2|\times T[0,1]\;.
\end{equation}
\end{remark}
$\,$\\\\A Poisson structure defines a 2-cocycle in the Chevalley-Eilenberg complex of $\textup{T}^*M.$ This integrates to a 2-cocycle on the 2-groupoid,\footnote{One should interpret this as the inverse of a higher groupoid van Est map.} which defines a multiplicative line bundle $\mathcal{L}_2\to \Pi_{2}(\textup{T}^*M):$
\begin{definition}\label{mult}
We have a 2-cocycle on $\Pi_2(\textup{T}^*M)$ given by 
\begin{equation}
    S:\Pi_{2}^{(2)}(\textup{T}^*M)\to\mathbb{R}\,,\;\; S[X]:=\int_{|\Delta^2|}X^*\Pi\;.
\end{equation}
We then have a multiplicative line bundle 
\begin{equation}
    \mathcal{L}_2\to \Pi_{2}(\textup{T}^*M)\,,\;\;\mathcal{L}_2=\Pi_{2}^{(1)}(\textup{T}^*M)\times\mathbb{C}\;,
\end{equation}
whose multiplication is defined as follows: for $\gamma_0,\gamma_1$ composing to $\gamma_2$ via $X,$ we define
\begin{equation}
    (\gamma_0,\lambda_0)\cdot_X(\gamma_1,\lambda_1)=(\gamma_2,\lambda\beta \,e^{\frac{i}{\hbar}\int_{|\Delta^2|}X^*\Pi}\,)\;.
\end{equation}
Sections of $\mathcal{L}_2\to  \Pi_{2}(\textup{T}^*M)$ are denoted $\Gamma(\mathcal{L}_2).$
\end{definition}
$\,$\\In the above definition, if $\int_{S^2}X^*\Pi\in 2\pi\hbar\mathbb{Z}$ for all morphisms $X:\textup{T}S^2\to T^*M,$ then the multiplication only depends on $\gamma_2.$ That we have defined a cocycle follows from Stokes' theorem.
\\\\Next, we define the twisted convolution algebra:
\begin{definition}
Consider two sections $w_1,\,w_2\in \Gamma(\mathcal{L}_2),$ ie. maps$\;\Pi_{2}^{(1)}(\textup{T}^*M)\to \mathbb{C}.$
Their twisted convolution
\begin{equation*}
w_1\ast w_2\in \Gamma(\mathcal{L}_2)
\end{equation*}
is given by\footnote{We are suppressing the dependence of $\ast$ on $\hbar$ in the notation.}
\begin{equation}
     (w_1\ast w_2)(\gamma)=\int_{\begin{subarray}{l}X\in \Pi^{(2)}_{2}(\textup{T}^*M)\\\gamma_2=\gamma\end{subarray}}w_1(\gamma_0)\cdot_X w_2(\gamma_1)\,\mathcal{D}X\;.
\end{equation}
 \end{definition}
Explicitly,
\begin{equation}
    (w_1\ast w_2)(\gamma)=\int_{\begin{subarray}{l}X\in \Pi^{(2)}_{2}(\textup{T}^*M)\\\gamma_2=\gamma\end{subarray}}w_1(\gamma_0)w_2(\gamma_1)e^{\frac{i}{\hbar}S[X]}\,\mathcal{D}X\;.
\end{equation}
This definition is formal and makes sense for any Lie 2-groupoid if we think of $X$ as just being a 2-morphism. If we apply it to a Lie 2-groupoid which is just a Lie 1-groupoid then (up to normalization) we recover the usual definition of the twisted convolution algebra (see eg. \cite{eli}).
\section{Main Result}
Here we discuss the three parts of \cref{formal}.
\subsection{Part 1: The 2-Groupoid Quantization Map}
\begin{definition}
We define the 2-groupoid quantization map by
\begin{equation}
    q:C^{\infty}(M)\to \Gamma(\mathcal{L}_2)\,,\;\;q_f(\gamma)=f(\gamma(1/2))
\end{equation}
(see \cref{mult} for definition of $\mathcal{L}_2$).
\end{definition}
The following  is part 1 of \cref{formal} and is a restatement of \cref{bon}, which we take as a fact for regular Poisson manifolds:
\begin{lemma}\label{eq}
Let $(M,\Pi)$ be a  Poisson manifold for which the right side of \ref{regstar} is smooth and let $f,g\in C_c^{\infty}(M).$ Then $(q_f\ast q_g)\vert_M=f\star g$ (here $\vert_M$ means the pullback to the space of objects).
\end{lemma}
\subsection{Part 2: The Geometric Quantization Map}
If for any morphism $X:\textup{T}S^2\to T^*M$ we have that
\begin{equation}\label{cond}
    \int_{S^2}X^*\Pi\in 2\pi\hbar\mathbb{Z}\;,
\end{equation}
then the multiplicative line bundle of \cref{mult} descends to a multiplicative line bundle over the source simply connected, symplectic groupoid $\Pi_1(\textup{T}^*M)$ (the underlying line bundle was constructed in \cite{bon}). 
\begin{definition}\label{multi1}
Recall that arrows in $\Pi_1(\textup{T}^*M)$ are equivalence classes of algebroid paths up to homotopy, relative to the endpoints (\cite{ruif}). We denote such an arrow by $[\gamma].$ Assuming \cref{cond}, the following defines a multiplicative line bundle $\mathcal{L}_1\to \Pi_1(\textup{T}^*M):$
\\\\The vector space over an arrow $[\gamma]$ 
is given by equivalence classes of points $(\gamma,\lambda)\in \mathcal{L}_2$ (defined in \cref{mult}), such that if $\gamma\sim \gamma',$ then
\begin{equation}
    (\gamma,\lambda)\sim (\gamma',\lambda\,e^{\frac{i}{\hbar}\int_D X^*\Pi}\,)\;,
\end{equation}
where $X:TD\to T^*M$ is any algebroid homotopy between $\gamma,\gamma'.$ The multiplication of $\mathcal{L}_2$ descends to a multiplication of $\mathcal{L}_1.$
\end{definition}
We now define the quantization map:
\begin{definition}\label{qant}
Assuming \cref{cond}, we define the geometric quantization map by
\begin{equation}
    Q:C_c^{\infty}(M)\to \Gamma(\mathcal{L}_1)\;,\;\;Q_f([\gamma])=\big(\gamma,\,\int_{X} q_f(\gamma')e^{\frac{i}{\hbar}\int_D X^*\Pi}\,\mathcal{D}X\,\big)\;,
\end{equation}
where the integral is over disks $X:TD\to TM$ with boundary components $\gamma,\gamma'$ ($\gamma$ is fixed). This definition is independent of the representative $\gamma$ of $[\gamma].$
\end{definition}
This definition explains the presence of Fourier-like transforms in the known cases where geometric quantization works, eg. our main example and \cite{pierre}, \cite{bone}, \cite{reiffelact}.
\subsection{Part 3: Representation on Sections of Prequantum Line Bundle}
On a pointed, simply connected symplectic manifold $(M,\omega, *),$ there is a canonical action of the multiplicative prequantum line bundle of $\mathcal{L}_2\to \Pi_2(\textup{T}^*M)$ on the prequantum line bundle (with connection) $\mathcal{L}_0\to M.$
\\\\The construction of the prequantum line bundle $\mathcal{L}_0$ is obtained by constructing the line bundle over $\Pi_1(T^*M)\cong \textup{Pair}\,M$\footnote{$\textup{Pair}\,M$ is the pair groupoid, it has a unique arrow between any two points in $M.$} as in the previous section, and pulling it back to the source fiber over $*.$ We will recall the standard construction in geometric quantization below (\cite{woodhouse}).
\\\\Assume that $[\omega]\in H^2(M,2\pi\hbar\mathbb{Z}).$ Consider the space of paths 
\begin{equation}
      P_0M=\{\gamma:[0,1]\to M\,:\,\gamma(0)=*\}\;.
    \end{equation}
Since $M$ is simpy connected, the map 
\begin{equation}
P_{0}M\to M\,,\;\;\;\gamma\mapsto\gamma(1)
\end{equation}
descends to a diffeomorphism on homotopy classes of paths relative to the endpoints, $P_{0}M/\sim.$ The trivial line bundle $P_0M\times\mathbb{C}$ descends to $M,$ where we quotient out by
\begin{equation}
    (\gamma_1,\lambda)\sim (\gamma_2,\lambda\,e^{\frac{i}{\hbar}\int_D\omega})\;,
\end{equation}
if $D$ is a disk whose boundary has components $\gamma_1, \gamma_2.$ We may therefore denote points in $M$ by $[\gamma]$ and points in the line bundle over $M$ by $[(\gamma,z)].$ 
\begin{definition}\label{multi0}
The prequantum line bundle $\mathcal{L}_0\to M$ is the line bundle defined by the previous construction. The trivial connection on $P_0M\times\mathbb{C}$ descends to a connection on $\mathcal{L}_0:$ let $\gamma_1,\gamma_2\in P_0M$ and let 
\begin{equation}
    \gamma:[0,1]\to M
\end{equation} 
be a path such that 
\begin{equation}
    \gamma(0)=\gamma_1(1)\,,\;\;\gamma(1)=\gamma_2(1)\;.
    \end{equation}
Parallel transport $\gamma_1(1)\to\gamma_2(1)$ over $\gamma$ sends 
\begin{equation}
[(\gamma_1,\lambda)]\to [(\gamma_2,\lambda\,e^{\frac{i}{\hbar}\int_D\omega})]\;,
\end{equation}
where $D$ is any disk whose boundary has components $\gamma_1,\gamma_2,\gamma.$ The space of sections is denoted $\Gamma(\mathcal{L}_0).$
\end{definition}
We are now ready to define the action of $\mathcal{L}_1$ on $\mathcal{L}_0:$
\begin{definition}\label{action}
We have an action of $\mathcal{L}_2$ on $\mathcal{L}_0,$ where
\begin{equation}
    (\gamma,\lambda )\in \mathcal{L}_2 =\Pi_2(\textup{T}^*M)\times\mathbb{C}
    \end{equation}
acts by parallel transport over $\gamma$\footnote{We are implicitly using the isomorphism $\omega:TM\to \textup{T}^*M,$ so that these algebroid morphisms into $T^*M$ are the same as maps of manifolds into $M.$} followed by multiplication by $\lambda.$ This descends to an action $\cdot$ of $\mathcal{L}_1\to \textup{Pair}\,M$ on $\mathcal{L}_0\to M.$ 
\end{definition}
Finally, we can discuss the representation:
\begin{definition}\label{reppp}
Let $2n=\textup{dim } M.$ We have a representation of $\Gamma(\mathcal{L}_1)$ on $\Gamma(\mathcal{L}_0)$ defined by
\begin{equation}
    (v\Psi)(m)=\int_{\{[\gamma]\in \textup{Pair}\,M:\gamma(1)=m\}} v([\gamma])\cdot\Psi(\gamma(0))\,\omega^n\;,
\end{equation}
where $v\in\Gamma(\mathcal{L}_1),\,\Psi\in\Gamma(\mathcal{L}_0).$
\end{definition}
\section{Example: Conventional Quantum Mechanics}\label{exam}
With the constructions made in the previous section, we are going to explicitly show what \cref{formal} gives in the example of conventional quantum mechanics. We will start with the K\"{a}hler structure and determine, in order: the quantization map, the non-perturbative Moyal product and the standard irreducible representations. The path integrals can be computed using the framework developed in \cite{Lackman4}, and we will assume the results (see appendix \ref{appen} for the idea).
\\\\Consider the pointed K\"{a}hler manifold $(M,\omega,J,*)=(\mathbb{R}^2,dp\wedge dq, p+iq,(0,0)).$\footnote{All of our constructions will depend only on this pointed K\"{a}hler structure.}
\begin{enumerate}
    \item  The symplectic groupoid is
\begin{equation}
  \Pi_1(\textup{T}^*\mathbb{R})=(\textup{Pair}\,\mathbb{R}^2\rightrightarrows\mathbb{R}^2,dp_1\wedge dq_1-dp_0\wedge dq_0)\;.
\end{equation}
We are identifying the fundamental groupoid of a simply connected space with its pair groupoid (which has a unique arrow between any two points).
\item The multiplicative line bundle $\mathcal{L}_1\to \textup{Pair}\,\mathbb{R}^2$ in \cref{multi1} is isomorphic to the trivial complex line bundle $\textup{Pair}\,\mathbb{R}^2\times\mathbb{C},$ twisted by the 2-cocycle given by
\begin{equation}
    \Omega(q,p,q_1,p_1,q_2,p_2)=\frac{1}{2}\begin{vmatrix} p_2-p & p_{1}-p   \\ q_{2}-q & q_{1}-q \end{vmatrix}\;,
\end{equation}
which is the (signed) area of the geodesic triangle with vertices $(q,p), (q_1,p_1), (q_2,p_2).$ This is the natural $2$-cocycle obtained using the K\"{a}hler structure.\footnote{The van Est map applied to $\Omega$ gives $\omega.$}
The multiplication of the line bundle is given by
\begin{equation}
(q_0,p_0,q_1,p_1,a)\cdot(q_1,p_1,q_2,p_2,b)=(q_0,p_0,q_2,p_2,ab\,\Omega(q_0,p_0,q_1,p_1,q_2,p_2))\;.
\end{equation}
\item There is a natural non-vanishing section of the prequantum line bundle $\mathcal{L}_0\to \textup{T}^*\mathbb{R}$ in \cref{multi0} obtained using the K\"{a}hler structure. It is given by
\begin{equation}
    (q,p) \mapsto (q,p,1)\;,
\end{equation}
where we identify $(q,p)$ with the geodesic connecting $0$ to $(q,p).$
This determines a trivialization of $\mathcal{L}_0,$ which we will assume from now on. In this trivialization, the canonical connection is given by 
\begin{equation}
    \frac{1}{2}(p\,dq-q\,dp)\;.
    \end{equation}
We have implicitly used the ``geodesic" trivialization of $\mathcal{L}_1$ in step 2.
\item The action of $\mathcal{L}_1$ on $\mathcal{L}_0$ in \cref{action}, obtained using parallel transport, is given by
\begin{equation}
    (q_0,p_0,p_1,q_1,a)\cdot(q_0,p_0,b)=(q_1,p_1,ab\,e^{\frac{i}{2\hbar}(p_1q_0-q_1p_0)})\;.
\end{equation}
We have that $(p_1q_0-q_1p_0)/2=\Omega(q_1,p_1,0,0,q_0,p_0).$
\item The representation of the corresponding twisted convolution algebra in \cref{reppp} on $L^2(\mathbb{R}^2)$ is given by
\begin{equation}\label{conv}
    (v\Psi)(q,p)=\frac{1}{4\pi\hbar}\int_{\mathbb{R}^2}v(q_0,p_0,q,p)\cdot \Psi(q_0,p_0)\,dp_0dq_0\;.
\end{equation}
\item The quantization map in \cref{qant} is given by\footnote{We have normalized by dividing out by \cref{kernel}.}
\begin{equation}\label{geoq}
   Q_f(q_0,p_0,q_1,p_1)=\frac{1}{4\pi\hbar}\int_{\mathbb{R}^2}f(q,p)e^{\frac{i}{2\hbar}[p(q_1-q_0)-q(p_1-p_0)]}\,dp\,dq\;.
\end{equation}
Geometrically, the exponent is equal to the (signed) area of the polygon determined by the points
\begin{equation}
   (0,0), (q_1,p_1), (q,p), (q_0,p_0)\;.
    \end{equation}
Thus, $Q_f$ equivariant under automorphisms (of the pointed K\"{a}hler structure).
\item Note that there is a Lagrangian polarization induced by the morphism
\begin{equation}
    \textup{Pair}\,\mathbb{R}^2\to \mathbb{R}^2\,,\;\;(q_0,p_0,q_1,p_1)\mapsto (q_1-q_0,p_1-p_0)
\end{equation}and $Q_f$ is polarized with respect to the connection $p_1d(q_1-q_0)-q_0d(p_1-p_0).$\footnote{Compare this with the construction of the quantization in \cite{eli}, for which a polarization needs to be chosen. The two quantization maps are related by an automorphism of $\mathbb{R}^2,$ $(x,y)\mapsto (-y,x).$}
\item The product $f\star g$ (defined in \cref{bon}) is given by
\begin{align}\label{my}
& (f\star g)(q,p)
\\&=\frac{1}{(4\pi\hbar)^2}\int_{\mathbb{R}^4}f(q_1,p_1)g(q_2,p_2)e^{\frac{i}{2\hbar}[(p_2-p)(q_1-q)-(q_2-q)(p_1-p)]}\,dp_1\,dq_1\,dp_2\,dq_2\;. 
\end{align}
This is equal to $(q_f\ast q_g)(q,p)$ in \cref{formal} and is the non-perturbative Moyal product (\cite{reiffelact}, \cite{Zachos}). It is invariant under automorphisms of the K\"{a}hler manifold and the integral kernel is $e^{\frac{i}{\hbar}\Omega},$ which is equal the integral kernel in \cref{kernel}.
\item We have that $Q_{f\star g}=\, Q_f\ast\, Q_g,$ and in fact 
\begin{equation}\label{equa}
    Q_f\Psi=f\star \Psi\;.
    \end{equation}
\\\textbf{This completes the constructions of \cref{formal}. We will now discuss irreducible subrepresentations of our $C^*$-algebra, using \cref{equa}.}
\item Covariantly constant Lagrangian polarizations are of the form 
\begin{equation}
    d\frac{\partial}{\partial q}-c\frac{\partial}{\partial p}\;,
\end{equation}
for constants $c,d\in\mathbb{C}$ (not both $0$). A polarized section can be written as
 \vspace{0.05cm}\\
    \begin{equation}\label{repp}
        \Psi_{c,d}(q,p)= e^{\frac{i}{2\hbar}(aq+bp)(cq+dp)}\psi(cq+dp)\;,
    \end{equation}
    \vspace{0.05cm}\\
for any $a,b$ such that 
\begin{equation}
\begin{vmatrix} a & c   \\ b & d \end{vmatrix}=1\;.\footnote{Note that $\{aq+bp,cq+dp\}=1,$ ie. these are canonical coordinates.}
\end{equation}
With respect to the K\"{a}hler polarization, polarized sections are of the form $e^{-\frac{|p+iq|^2}{2\hbar}}\psi(p+iq)\,.$ 
\\\\\textbf{These give irreducible subrepresentations of the left representation of $\star,$ ie.} $Q_f\Psi_{c,d}$ \textbf{is of the form \ref{repp}, for the same $c,d$ and for any $f.$ See \cref{ch}.}
\item The actions of the quantizations of $aq+bp,\,cq+dp$ on $\Psi_{c,d}$ (as in \ref{repp}) are given by\footnote{Our conventions imply that $[Q_q,Q_p]= 2i\hbar .$}
\begin{align}
    & Q_{aq+bp}\Big[e^{\frac{i}{2\hbar}(aq+bp)(cq+dp)}\psi(cq+qp)\Big]=e^{\frac{i}{2\hbar}(aq+bp)(cq+dp)}\Big[\frac{\hbar}{i}\psi'(cq+dp)\Big]\,,
    \\& Q_{cq+dp} \Big[e^{\frac{i}{2\hbar}(aq+bp)(cq+dp)}\psi(cq+qp)\Big]=e^{\frac{i}{2\hbar}(aq+bp)(cq+dp)}\big[(cq+dp)\psi(cq+dp)\big]\,.
\end{align}
We can read off the position, momentum and Segal-Bargmann representations (ie. the holomorphic representation \cite{hall}).\footnote{The holomorphic representation is the only one which forms a $C^*$-algebra subrepresentation, the former are only algebra subrepresentations, but there is essentially a unique measure turning them into $C^*$-algebra representations.} 
\\\\We are describing the representation of the quantum observables on all vector states\footnote{A pure state is one which generates an irreducible representation, these polarized sections are pure states. However, the quantization map is defined on all sections of the prequantum line bundle, and some of these sections are mixed states. in the sense of $C^*$-algebras. See \cite{Gleason}.} before choosing a polarization. In traditional quantization, one chooses a polarization first and the Kostant-Souriau quantization prescription only quantizes a small subspace of observables which preserves polarized sections. 
\\\\Generally, if $\Psi$ is a polarized section, then
\begin{equation}\label{ch}
Q_f\Psi(q,p)=\frac{1}{4\pi\hbar}\int_{\mathbb{R}^2}f(q',p')\Psi(q'-q,p'-p)e^{\frac{i}{2\hbar}(qp'-pq')}\,dp'dq'\;.
\end{equation}
\item The involution $q+ip\mapsto q-ip$ 
lifts to an isomorphism of representations from states polarized with respect to $\partial_p$ to states polarized with respect to $\partial _q.$ The isomorphism is given simply by integration with respect to $dq,$ where we identify the vector spaces over $(q',p)$ and $(q,p)$ using the flat connection along the polarization, ie. under parallel transport over $(q',p)\to (q,p),$ we have that
\begin{equation}
    (q',p,a)\mapsto (q,p,ae^{-\frac{i}{2\hbar}p(q'-q)})\;.
\end{equation}
Performing the integral, we get
\begin{equation}
    \Big((q,p)\mapsto e^{-\frac{i}{2\hbar}pq}\psi(q)\Big)\mapsto \Big((q,p)\mapsto e^{\frac{i}{2\hbar}qp}\int_{\mathbb{R}}e^{-\frac{i}{\hbar}pq'}\psi(q')\,dq'\;\Big)\;.
\end{equation}
\end{enumerate}
\subsection{The Classical Limit}
We can analyze the classical limit from the 2-groupoid perspective. First, note that given a Poisson structure $\Pi$ on $M,$ we get a natural Poisson structure on $M\times \mathbb{R},$ where the induced Poisson structure on $M\times \{\hbar\}$ is $\hbar\Pi$.\footnote{We can replace $\mathbb{R}$ with any subset with $0$ as an accumulation point.} Then, what we have really done in this paper is integrated the extension of Lie algebroids\footnote{$  \mathbb{C}_{M\times \mathbb{R}}$ is the trivial line bundle over $M\times \mathbb{R}.$}
\begin{equation}
   \mathbb{C}_{M\times \mathbb{R}}\to T^*(M\times \mathbb{R})\times\mathbb{C}\to T^*(M\times\mathbb{R})
\end{equation}
to a multiplicative line bundle
\begin{equation}
    \mathcal{L}\to \Pi_2(T^*(M\times\mathbb{R}))\;.
\end{equation}This extension is the one obtained by interpreting $\Pi$ as a 2-cocycle on $T^*(M\times\mathbb{R}).$
\\\\The $\hbar=0$ case is special: isomorphism classes (ie. orbits) of objects in the 2-groupoid of $\hbar\Pi$ are its symplectic leaves, and $\hbar=0$ is the unique case where the symplectic leaves are just the points in $M.$ As a result, the pure states at $\hbar=0$ are just pointwise evaluation. From the 2-groupoid perspective, the measurement problem is due to the existence of isomorphisms between distinct points in phase space of a symplectic manifold, for $\hbar\ne 0.$
\\\\Even though algebroid paths are constant at $\hbar=0,$ there is still nontrivial dynamics since we aren't varying the 2-cocycle $\Pi$ with $\hbar,$ thus the $\hbar=0$ case still remembers the Poisson structure. This is the only case where the $2$-cocycle $\Pi$ isn't proportioinal to the canonical 2-cocycle on $\textup{T}^*M$ (since when $\hbar=0,$ the canonical 2-cocycle is $0$).
\\\\One can take a line bundle over $M$ which prequantizes each of its symplectic leaves. A representation of the quantization on this space of sections will pull back to a representation on sections of the line bundle over each symplectic leaf (a ``subrepresentation"). It is natural to take a subspace of sections over $M$ which pull back to irreducible representations over each leaf (eg. by taking a foliation by coisotropic submanifolds whose intersection with each symplectic leaf is Lagrangian). We can think of  this as a classical mixture of quantum systems.
\section{Concluding Remarks}
In our main example, one can check that if $f$ is periodic and $\Psi$ transforms as a theta function under translations, then $f\star \Psi$ transforms as a theta function. Therefore, if we carry out our construction on the symplectic torus, we will get a representation of the noncommuative torus on sections of the prequantum line bundle. It would be interesting to construct the irreducible representations.
\\\\In \cite{schl} (see also \cite{bord}, \cite{chan}) it is shown that there is a sequence of ``quantization maps" and representations associated to a compact K\"{a}hler manifold (Berezin-Toeplitz quantization) which induce a formal deformation quantization. However, these quantization maps are far from injective. We expect that our construction can be carried out in this context, and that their representations factor through our geometric quantization map $Q,$ with the path integrals being computed as in \cite{Lackman4}.
\begin{appendices}\label{appen}
We will briefly describe the definition of path integrals over spaces of Lie algebroid morphisms given in \cite{Lackman4}, specialized to the example of the star product in the form \cref{regstar}, ie.
\begin{equation*}
    (f\star g)(m)=\int_{\pi(X(\infty))=m}f(\pi(X(1)))\,g(\pi(X(0)))\,e^{\frac{i}{\hbar}\int_{D}X^*\Pi}\,\mathcal{D}X\;.
    \end{equation*}
The steps are:
\begin{enumerate}
    \item triangulate the disk (with the marked points contained in the vertices) and form the associated simplicial set $\Delta_D,$
    \item integrate the Poisson manifold to the local symplectic groupoid,
    \item approximate the space of Lie algebroid morphisms by morphisms of simplicial sets between $\Delta_D$ and the local symplectic groupoid,
    \item integrate the Poisson structure to a 2-cocycle, via the van Est map,
    \item form the approximations of the action (generalized Riemann sums),
    \item construct the measure on the hom space of simplicial sets by using the symplectic form on the symplectic leaves and the Haar measures on the isotropy groups,
    \item compute the approximations of the functional integral,
    \item take a limit over triangulations of the disk.
\end{enumerate}   
The path integrals in \ref{exam} can then be computed using the integral representation of the delta function.
\end{appendices}


\begin{thebibliography}{9}
\bibitem{pierre}
Pierre Bieliavsky, Philippe Bonneau, Yoshiaki Maeda. \textit{Universal Deformation Formulae, Symplectic Lie groups and Symmetric Spaces.} Pacific Journal of Mathematics. Volume 230 No. 2 April 2007.
 \bibitem{bon}
F. Bonechi, A. S. Cattaneo and M. Zabzine, \textit{Geometric quantization and
non-perturbative Poisson sigma model.} Adv. Theor. Math. Phys. 10 (2006)
683 [arXiv:math/0507223].
\bibitem{bone}
F. Bonechi, N. Ciccoli, N. Staffolani, M. Tarlini. \textit{The quantization of the symplectic groupoid of the standard Podle\`{s} sphere.} Journal of Geometry and Physics, (2010), 62, 1851-1865.
\bibitem{bord}
M. Bordeman, E. Meinrenken and M. Schlichenmaier. \textit{Toeplitz quantization of K\"{a}hler 
manifolds and gl(n), n $\to\infty$ 
limits.} Comm. Math. Phys. 165 (1994), 281-296. 
 \bibitem{catt}
 A. S. Cattaneo and G. Felder. \textit{A Path Integral Approach to the Kontsevich Quantization Formula}. Comm Math Phys 212, 591–611 (2000). https://doi.org/10.1007/s002200000229
 \bibitem{chan}
K. Chan, N.C. Leung and Q. Li. \textit{Quantization of Kähler manifolds.} (2020) Journal of Geometry and Physics, 163, 104143.
\bibitem{ruif}
Marius Crainic and Rui Loja Fernandes. \textit{Integrability of Poisson Brackets.} J. Differential Geom. 66 (1) 71 - 137, January, 2004. https://doi.org/10.4310/jdg/1090415030
\bibitem{fedosov}
Boris V. Fedosov. \textit{A simple geometrical construction of deformation quantization.} J. Differential Geom. 40(2): 213-238 (1994). DOI: 10.4310/jdg/1214455536
\bibitem{Gleason}
Jonathan Gleason. \textit{The C*-algebraic formulation of quantum mechanics}. (2009) \href{https://math.uchicago.edu/~may/VIGRE/VIGRE2009/REUPapers/Gleason.pdf}{THE C*-ALGEBRAIC FORMALISM OF QUANTUM MECHANICS}
\bibitem{guk}
Sergei Gukov and Edward Witten. \textit{Branes and quantization.}
Adv. Theor. Math. Phys. 13(5): 1445-1518 (2009).
 \bibitem{eli}
Eli Hawkins. \textit{A Groupoid Approach to Quantization.}
J. Symplectic Geom. 6 (2008), no. 1, 61-125.
\bibitem{eli2}
Eli Hawkins. \textit{An Obstruction to Quantization of the Sphere.} Communications in Mathematical Physics 283 (2007): 675-699.
\bibitem{gamme}
 A. Gammella. \textit{Tangential products.} Letters in Mathematical Physics 51 , 1–15 (2000).
\bibitem{hall}
B. Hall. \textit{Geometric Quantization and the Generalized Segal--Bargmann Transform for Lie Groups of Compact Type.} Commun. Math. Phys. 226, 233–268 (2002). 
\bibitem{kontsevich}
M. Kontsevich. \textit{Deformation Quantization of Poisson Manifolds.} Letters in Mathematical Physics 66, 157–216 (2003). https://doi.org/10.1023/B:MATH.0000027508.00421.bf
\bibitem{Lackman3}
Joshua Lackman. \textit{A Formal Equivalence of Deformation Quantization and Geometric Quantization (of Higher Groupoids) and Non-Perturbative Sigma Models.} \href{https://arxiv.org/abs/2303.05494}{arXiv:2303.05494} (2023).
\bibitem{Lackman4}
Joshua Lackman. \textit{A Groupoid Construction of Functional Integrals: Brownian Motion and Some TQFTs.}  \href{https://arxiv.org/abs/2402.05866}{arXiv:2402.05866v2} [math.DG] (2024).
\bibitem{nolle}
C. N\"{o}lle. \textit{On the relation between geometric and deformation quantization.} \href{https://arxiv.org/abs/0809.1946}{arXiv:0809.1946} (2008).
\bibitem{nolle2}
C. N\"{o}lle. \textit{Geometric and deformation quantization.} \href{https://arxiv.org/abs/0903.5336}{arXiv:0903.5336v2} [math-ph] (2009).
\bibitem{rieffel}
M. Rieffel. \textit{Questions on quantization.} Operator Algebras and Operator Theory (1997, Shanghai), Contemp. Math. 228 (1998), 315–326, quant-ph/9712009.
\bibitem{rieffel2}
M. Rieffel. \textit{Lie group convolution algebras as deformation quantizations of linear Poisson
structures.} Amer. J. Math. 112 (1990), 657–685. MR1064995.
\bibitem{reiffelact}
M. Rieffel. \textit{Deformation Quantization for Actions of $R^d$}. Memoirs of the American Mathematical Society, Volume: 106; 1993; 93 pp MSC: Primary 46; Secondary 35.
\bibitem{schl}
M. Schlichenmaier. \textit{Deformation quantization of compact Kähler manifolds by Berezin-Toeplitz quantization.} In: Dito, G., Sternheimer, D. (eds) Conf\'{e}rence Mosh\'{e} Flato 1999. Mathematical Physics Studies, vol 21/22. Springer, Dordrecht. 
\bibitem{woodhouse}
D.J. Simms, N.M.J. Woodouse. \textit{Lectures on Geometric Quantization.} 1976, Volume 53
ISBN : 978-3-540-07860-9
\bibitem{wein}
A. Weinstein. \textit{Tangential deformation quantization and polarized
symplectic groupoids.} Deformation Theory and Symplectic
Geometry, S. Gutt, J. Rawnsley, and D. Sternheimer, eds.,
Mathematical Physics Studies 20, Kluwer, Dordrecht, 1997,
301-314.
\bibitem{weinstein}
Alan Weinstein. \textit{Symplectic Groupoids, Geometric Quantization, and Irrational Rotation Algebras.} In: Dazord, P., Weinstein, A. (eds) Symplectic Geometry, Groupoids, and Integrable Systems. Mathematical Sciences Research Institute Publications, vol 20,  (1991). Springer, New York, NY.
\bibitem{zhuc}
Chenchang Zhu. \textit{Lie II theorem for Lie algebroids via higher Lie groupoids.} (2010)
\href{https://arXiv.org/abs/math/0701024v2}{arXiv:math/0701024v2 [math.DG]}
\bibitem{Zachos}
Cosmas Zachos. \textit{Geometrical evaluation of star products}. J. Math. Phys. 41, 5129-5134 (2000) https://doi.org/10.1063/1.533395.
\end{thebibliography}
 \end{document}